# AVERAGE-CASE COMPLEXITY AND DECISION PROBLEMS IN GROUP THEORY

ILYA KAPOVICH, ALEXEI MYASNIKOV, PAUL SCHUPP, AND VLADIMIR SHPILRAIN

ABSTRACT. We investigate the average-case complexity of decision problems for finitely generated groups, in particular the word and membership problems. Using our recent results on "generic-case complexity" we show that if a finitely generated group $G$ has the word problem solvable in subexponential time and has a subgroup of finite index which possesses a non-elementary word-hyperbolic quotient group, then the average-case complexity of the word problem for $G$ is linear time, uniformly with respect to the collection of all length-invariant measures on $G$. For example, the result applies to all braid groups $B_n$.

## 1. INTRODUCTION

One of the most classical decision problems is the word problem for groups. If $G$ is a group with a finite generating set $A = \{x_1, \ldots, x_k\}$, then the *word problem $WP(G, A)$ of $G$ relative to $A$* consists of all those words $w$ in the alphabet $X = A \cup A^{-1}$ which represent $1 \in G$. Thus $WP(G, A) \subseteq X^*$, where $X^*$ is the set of all words in the alphabet $X$. For finitely generated groups the worst-case complexity of the word problem does not depend on the choice of a finite generating set and for that reason the reference to a generating set is often suppressed.

In general, group word problems may be very complicated with respect to worst-case complexity. The basic result of Novikov and of Boone (see for example [35, 36, 1]) states that there exists a finitely presented group $G$ such that for any finite generating set $A$ of $G$ the word problem $WP(G, A)$ is undecidable. However, "most" finitely presented groups are word-hyperbolic and hence have the word problem solvable in linear time. This fact was first observed by Gromov [27] and made precise by Ol'shanskii [38] and Champetier [14]. They introduced the notion of a "generic" group-theoretic property and that notion is being successfully explored by many authors [4, 5, 6, 7, 15, 16, 17, 18, 28, 29].

In our recent paper [32] we introduced the notion of *generic-case complexity* and showed that for most of the groups usually studied in combinatorial group theory, the generic-case complexity of the word problem is linear time, meaning that there is a partial algorithm which gives the correct answer for "most" inputs of the problem in linear time. We will recall the precise definition of generic-case complexity later in Section 2. Such a result can hold for groups where the worst-case complexity of the word problem is very high or the problem is even undecidable because this notion completely disregards the behavior of the algorithm on a "small" set of "difficult" inputs.

In the present article we show that our results regarding generic-case complexity can in fact be used to obtain precise average-case complexity statements, which concern the expected value of the complexity over the whole set of inputs, including the "difficult" ones. The basic idea is very straightforward and is often used in practice. If we have a total algorithm $M_1$ solving a decision problem $\mathcal{D}$ whose worst-case complexity is not "too high" and we also have a partial algorithm $M_2$ solving the problem with "strongly" low generic-case complexity, then by running $M_1$ and $M_2$ in parallel we have a total algorithm $M_1 || M_2$ for which we can hope to prove low average-case complexity.

There are few average-case complexity results about decision problems related to finitely presented groups. We can mention here two papers by Wang [41, 43], whose excellent results are very different from ours both in substance and in the technique used (for example, in Wang's approach the set of instances of a problem involves all finite group presentations, rather than a fixed one).

**Definition 1.1.** Let $X$ be a finite alphabet with $|X| \geq 2$ elements. A *discrete probability measure* on $X^*$ is a function $\mu : X^* \to [0,1]$ such that $\sum_{w \in X^*} \mu(w) = 1$.

We will say that $\mu$ is *length-invariant* if for any words $w, w' \in X^*$ with $|w| = |w'|$ we have $\mu(w) = \mu(w')$.

Requiring that a measure be length-invariant is a very natural assumption, since most complexity classes are defined in terms of the length of an input word.

**Convention 1.2.** We follow [37] for our conventions on computational complexity. In particular, by an *algorithm* we will always mean a deterministic multi-tape Turing machine. In this paper we restrict our consideration to deterministic time-complexity classes, although the methods of this paper can be applied to the analysis of average-case behavior of more general complexity classes. We also assume that for a time-complexity class $\mathcal{C}$ the collection of functions bounding the time of a computation consists of proper complexity functions $f(n) \geq n$ and that for any function $f(n)$ in this collection and for any integer $C \geq 1$ the function $Cf(Cn+C)+C$ also belongs to this collection.

**Definition 1.3** (Subexponential functions)**.** We say that a non-negative function $f(n)$ is *subexponential* if for any $r > 1$ we have

$$\lim_{n \to \infty} \frac{f(n)}{r^n} = 0.$$

Note that this implies that for every $r > 1$

$$\sum_{n=1}^{\infty} \frac{f(n)}{r^n} < \infty.$$

**Definition 1.4** (Average-case complexity). Let $X$ be a finite alphabet with $|X| \geq 2$ elements. Let $\mathcal{D} \subseteq X^*$ be a language and let $M$ be an algorithm which for every $w \in X^*$ decides whether or not $w \in \mathcal{D}$ in time $T(w) < \infty$. Let $f(n)$ be a non-decreasing positive function. Let $\mu : X^* \to [0,1]$ be a discrete probability measure.

(1) We will say that $M$ *solves $\mathcal{D}$ with average case time-complexity bounded by $f(n)$ relative to $\mu$* if

$$\int_{X^*} \frac{T(w)}{f(|w|)} \mu(w) = \sum_{w \in X^*} \frac{T(w)}{f(|w|)} \mu(w) < \infty.$$

If $f(n)$ satisfies the bound constraint of a time-complexity class $\mathcal{C}$ we will say that $M$ *solves $\mathcal{D}$ with average case time-complexity in $\mathcal{C}$ relative to $\mu$*.

(2) Let $\Re$ be a family of discrete probability measures on $X^*$. We say that $M$ *solves $\mathcal{D}$ with average case time-complexity bounded by $f(n)$ uniformly relative to $\Re$* if there is $0 < C < \infty$ such that for any $\mu \in \Re$

$$\int_{X^*} \frac{T(w)}{f(|w|)} \mu(w) = \sum_{w \in X^*} \frac{T(w)}{f(|w|)} \mu(w) \leq C.$$

If $f(n)$ satisfies the constraint of a time-complexity class $\mathcal{C}$ we will say that $M$ *solves $\mathcal{D}$ with average case time-complexity in $\mathcal{C}$ uniformly relative to $\Re$*.

**Remark 1.5.** Suppose $M$ solves $\mathcal{D}$ with average case time-complexity bounded by $f(n)$ relative to $\mu$. Then

$$\sum_{w \in X^*} \frac{T(w)}{f(|w|)} \mu(w) = \sum_{n=0}^{\infty} \sum_{w \in X^*, |w|=n} \frac{T(w)\mu(w)}{f(n)} < \infty.$$

Therefore

$$\sup_{n \geq 0} \frac{\sum_{w \in X^*, |w|=n} T(w)\mu(w)}{f(n)} < \infty,$$

and

$$\sum_{w \in X^*, |w|=n} T(w)\mu(w) = O(f(n)).$$

**Convention 1.6.** We shall denote by $\mathcal{S}ub\mathcal{E}xp$ the class of languages decidable in deterministic subexponential time.

We are now ready to state the main result of this paper. While the concept of a group being non-amenable plays an important role in our results, the reader need only have in mind that any group containing a free group of rank at least 2 is non-amenable.

**Theorem A.** *Let $G$ be a finitely presented group where the word problem is in $\mathcal{S}ubExp$. Suppose $G$ has a subgroup of finite index which possesses a non-amenable quotient group $\bar{G}_1$ whose word problem is solvable in a complexity class $\mathcal{C}$, where $\mathcal{C} \subseteq \mathcal{S}ubExp$.*

*Then the word problem $WP(G, A)$ for $G$ is solvable with average-case complexity in $\mathcal{C}$ uniformly relative to the family of all length-invariant discrete probability measures $\mu : (A \cup A^{-1})^* \to [0, 1]$.*

We can now formulate some more concrete corollaries of the above theorem.

**Corollary B.** *Let $G$ be a finitely presented group where the word problem is solvable in subexponential time. Let $A$ be a generating set of $G$ and let $X = A \cup A^{-1}$. Let $\Re$ be the family of all length-invariant discrete probability measures $\mu : X^* \to [0, 1]$.*

1. *Suppose that $G$ has a subgroup of finite index that possesses a non-elementary word-hyperbolic quotient group. Then the word problem $WP(G, A) \subseteq X^*$ is solvable with linear average-case time-complexity uniformly relative to $\Re$.*
2. *Suppose that $G$ has a subgroup of finite index that possesses a non-amenable automatic quotient group. Then the word problem*

$$WP(G, A) \subseteq X^*$$

   *is solvable with quadratic average-case time-complexity uniformly relative to $\Re$.*

From the point of view of group theory, one weakness of average-case complexity is that it says nothing about the complexity of the considered problem for finitely generated subgroups. If $G$ is a finitely presented group where the word problem is solvable in subexponential time, then an easy argument shows that the word problem is solvable in subexponential time for the direct product $G_1 = G \times F(a, b)$. Note that the free group $F(a, b)$ is non-elementary word-hyperbolic and $G_1$ admits a homomorphism onto $F(a, b)$. Thus Corollary B implies that any finitely presented group with word problem solvable in subexponential time can be embedded into a finitely presented group which has word problem solvable with linear-time average-case complexity (relative to the collection of all length-invariant probability measures).

**Example 1.7.** Let $B_n$ be the $n$-strand braid group. We observed in [32] that for any $n \geq 3$ the group $B_n$ has word-problem solvable with strongly linear-time generic-case complexity. The reason for this is that the pure-braid group $P_n$ (which has finite index in $B_n$) admits a homomorphism onto the group $P_3 \cong F(a, b) \times \mathbb{Z}$ (pulling out all but the first three strands of a pure braid), and the group $P_3$ in turn maps onto a non-elementary word-hyperbolic group $F(a,b)$. Since $B_n$ is automatic [21], the word problem in $B_n$ is solvable in quadratic time. Hence by Corollary B the word problem in $B_n$ is solvable with linear time average-case complexity (uniformly relative to the collection of all length-invariant measures for any fixed finite generating set of $B_n$).

**Example 1.8.** Let $G$ be a finitely generated linear group over a field of characteristic zero. Then the result of Lipton and Zalcstein [34] shows that the word problem in $G$ is solvable in log-space and thus in polynomial time. If $G$ has a finite index subgroup with a non-elementary word-hyperbolic quotient, then by Corollary B the word problem in $G$ is solvable with linear time average-case complexity.

**Example 1.9.** Suppose $G$ is a group which can be given by a finite presentation involving at least two more generators than defining relators. By the result of Baumslag and Pride [9] $G$ has a subgroup of finite index that maps homomorphically onto the free group of rank two. Suppose also that $G$ has the word problem solvable in subexponential time. Then by part 1 of Corollary B the word problem in $G$ is solvable in linear time on average.

Similar results hold for the *subgroup membership problem*.

Recall that if $G$ is a group with a finite generating set $A = \{x_1, \ldots, x_k\}$ and $H \leq G$ is a subgroup, then the *membership problem $MP(G, H, A)$ for $H$ in $G$ relative to $A$* consists of all those words $w$ in the alphabet $X = A \cup A^{-1}$ which represent elements of $H$. Thus $MP(G, H, A) \subseteq X^*$. As in the case of the word problem, the complexity of the membership problem does not depend on the choice of $A$.

Our main result regarding the membership problem is:

**Theorem C.** *Let $G$ be a finitely presented group and let $H \leq G$ be a subgroup where the membership problem for $H$ in $G$ is in $\mathcal{S}ubExp$. Let $G_1 \leq G$ be a subgroup of finite index in $G$ such that $H \leq G_1$ and let $\phi : G \to \bar{G}$ be an epimorphism.*

*Suppose there is a subgroup $\phi(H) \leq K \leq \bar{G}$ such that the Schreier coset graph for $\bar{G}$ over $K$ is non-amenable and such that the membership problem for $K$ in $\bar{G}$ is solvable in complexity class $\mathcal{C} \subseteq \mathcal{S}ubExp$.*

*Then for any finite generating set $A$ of $G$ the membership problem*

$$MP(G, H, A)$$

*for $H$ in $G$ is solvable with average-case complexity in $\mathcal{C}$ uniformly relative to the family of all length-invariant discrete probability measures $\mu : (A \cup A^{-1})^* \to [0, 1]$.*

**Corollary D.** *Let $G$ be a finitely presented group and let $H \leq G$ be a subgroup where the membership problem for $H$ in $G$ is solvable in subexponential* *time. Let $G_1 \le G$ be a subgroup of finite index in $G$ such that $H \le G_1$. Let $\phi : G_1 \to \bar G$ be an epimorphism and let $\phi(H) \le K \le \bar G$.*

1. *Suppose $\bar G$ is a non-elementary word-hyperbolic and that $K \le \bar G$ is a rational subgroup of infinite index.*

    *Then for any finite generating set $A$ of the group $G$ the membership problem $MP(G, H, A)$ for $H$ in $G$ is solvable with linear-time average-case complexity uniformly relative to the family of all length-invariant discrete probability measures $\mu : (A \cup A^{-1})^* \to [0, 1]$.*
2. *Suppose $\bar G$ is automatic and that $K \le \bar G$ is a rational subgroup such that the Schreier coset graph for $\bar G$ over $K$ is non-amenable (and hence $K$ has infinite index in $\bar G$). Then for any finite generating set $A$ of $G$ and for any length-invariant discrete probability measure $\mu : (A \cup A^{-1})^* \to [0, 1]$ the membership problem $MP(G, H, A)$ for $H$ in $G$ is solvable with quadratic-time average-case complexity uniformly relative to the family of all length-invariant discrete probability measures $\mu : (A \cup A^{-1})^* \to [0, 1]$.*

We should make an important disclaimer. The notion of average-case complexity used in this paper is rather rough and certainly does not have the many robustness properties that are desirable at finer levels of the complexity theory. For example, while our results are independent of the choice of a finite generating set for a group, more delicate issues such as, for instance, model independence still have to be addressed (see [33, 43] for a more detailed discussion).

Nevertheless, we believe that our results constitute a valuable step in studying the largely unexplored field of average-case complexity of group-theoretic algorithmic problems.

## 2. Generic-case complexity

We refer the reader to the classical text of Papadimitriou [37] for the basic definitions and conventions regarding computational complexity. We need to recall some definitions from our earlier paper [32].

**Definition 2.1** (Asymptotic density). Let $X$ be a finite alphabet with at least two elements and let $X^*$ denote the set of all words in $X$. Let $S$ be a subset of $X^*$. For every $n \ge 0$ let $B_n$ be the set all words in $A^*$ of length at most $n$.

The *asymptotic density* $\rho(S)$ for $S$ in $X^*$ is defined as

$$\rho(S) := \limsup_{n \to \infty} \rho_n(S)$$

where

$$\rho_n(S) := \frac{|S \cap B_n|}{|B_n|},$$

If the actual limit $\lim_{n \to \infty} \rho_n(S)$ exists, we write $\hat\rho(S) := \rho(S)$.

The above notion was first suggested in the paper of Borovik, Myasnikov and Shpilrain [12].

If $a_n \geq 0$ and $\lim_{n\to\infty} a_n = 0$, we will say that the convergence is *exponentially fast* if there is $0 \leq \sigma < 1$ and $C > 0$ such that for every $n \geq 1$ we have $a_n \leq C\sigma^n$. Similarly, if $\lim_{n\to\infty} b_n = 1$ (where $0 \leq b_n \leq 1$), we will say that the convergence is *exponentially fast* if the limit $\lim_{n\to\infty}(1 - b_n) = 0$ converges exponentially fast.

**Definition 2.2** (Generic performance of a partial algorithm). Let $X$ be a finite alphabet with at least two letters and let $\mathcal{D} \subseteq X^*$ be a language in $X$.

Let $\Omega$ be a correct partial algorithm for $\mathcal{D}$ which accepts as inputs words from $X^*$. That is, whenever $\Omega$ reaches a (positive or negative) decision on whether a word $w$ belongs to $\mathcal{D}$, that decision is correct. Let $\mathcal{C}$ be a complexity class (e.g. linear time, quadratic time, linear space, etc).

We say that $\Omega$ *solves $\mathcal{D}$ with generic-case complexity $\mathcal{C}$* if there is a subset $S \subseteq X^*$ with $\hat{\rho}(S) = 1$ such that for every $w \in S$ the partial algorithm $\Omega$ decides whether or not $w$ is an element of $\mathcal{D}$ within the complexity bound $\mathcal{C}$ (in terms of $w$). If in addition $\lim_{n\to\infty} \rho_n(S) = 1$ converges exponentially fast, we say that $\Omega$ *solves $\mathcal{D}$ with generic-case complexity strongly $\mathcal{C}$*.

Note that unlike the average-case complexity defined in the Introduction, generic-case complexity totally disregards the complexity of the algorithm on the "small" set of inputs $X^* - S$.

## 3. Proofs of the main results

**Convention 3.1.** Let $G$ be a group and let $A$ be a finite alphabet equipped with a map $\pi : A \to G$ such that $\pi(A)$ generates $G$. In this case we say that $(A, \pi)$ (or, by abuse of notation, just $A$) is a *finite generating set* for $G$. Thus we allow different letters of $A$ to represent the same element of $G$ and we also allow some letters of $A$ to represent $1 \in G$. Every word $w$ in the alphabet $A \cup A^{-1}$ represents an element of $G$ which we will still denote by $\pi(w)$.

If $H \leq G$ is a subgroup, we define the *Schreier coset graph* $\Gamma(G, H, A)$ as follows. The vertices of $\Gamma$ are cosets $\{Hg \,|\, g \in G\}$. The graph $\Gamma$ is oriented. The edge-set of $\Gamma$ is partitioned as $E\Gamma = E^+ \sqcup E^-$, where the set $E^+$ is in one-to-one correspondence with $V\Gamma \times A$. Namely for each coset $Hg$ and each $a \in A$ there is an associated edge $e$ from $Hg$ to $Hga$ with label $a$ in $E^+$. The set $E^-$ consists of formal inverses of the edges from $E^+$, where the inverse of an the edge $e$ above is the edge $e^{-1}$ from $Hga$ to $Hg$ with label $a^{-1}$. Thus the edges of $\Gamma$ are labeled by letters of $A \cup A^{-1}$. Note that if $k = |A|$ then $\Gamma$ is connected and $2k$-regular.

Intuitively, for a regular connected graph $\Gamma$ is to be *non-amenable* corresponds to the graph "growing fast in all directions". We refer the reader to [8, 13, 19, 26, 31, 40, 44, 45] for a detailed discussion of non-amenability for groups and graphs, including the many equivalent definitions (such as the Følner condition, the growth-rate criterion for regular graphs, the doubling condition etc). A finitely generated group $G$ is *non-amenable* if for some (and therefore for any) finite generating set $A$ of $G$ the Cayley graph $\Gamma(G, A)$ is amenable. This is equivalent to the standard definition of amenability for finitely generated groups. That is, a finitely generated group $G$ is amenable if and only if for any action of $G$ by homeomorphisms on a compact space $Q$ there exists a $G$-invariant probability measure on $Q$.

Recall that a subgroup $H$ of an automatic group $G$ is said to be *rational* if for some automatic language $L$ for $G$ the pre-image $L_H$ of $H$ in $L$ is itself a regular language. If $G$ is word-hyperbolic then rational subgroups are also often called *quasiconvex*. We refer the reader to [2, 10, 20, 22, 23, 24, 27, 25, 21, 39] for background information on hyperbolic and automatic groups and their rational subgroups. For the moment we need only recall that the word problem is solvable in quadratic time for any automatic group and in linear time for any hyperbolic group [2, 3, 21, 23, 30]. Also, if $H$ is a rational subgroup of a hyperbolic ( automatic) group $G$ then the membership problem for $H$ in $G$ is solvable in linear (quadratic) time.

Our main technical tool is:

**Proposition 3.2.** *Let $X$ be a finite alphabet with at least two elements and let $\mathcal{D} \subseteq X^*$. Suppose that $\mathcal{D}$ is decidable in time $f(n)$ and strongly generically decidable in time $f_1(n) \leq f(n)$ where the function $f(n)/f_1(n)$ is subexponential. Then the language $\mathcal{D}$ is decidable with average-case time-complexity bounded by $f_1(n)$ uniformly relative to the family of all length-invariant discrete probability measures $\mu : (A \cup A^{-1})^* \to [0, 1]$.*

*Proof.* Let $M'$ be an algorithm solving $\mathcal{D}$ in time $f(n)$. Let $M''$ be a partial algorithm which solves $\mathcal{D}$ strongly generically in time $f_1(n)$. Define $M$ to be the algorithm consisting of running $M'$ and $M''$ concurrently. Let $\mu : X^* \to [0, 1]$ be a length-invariant discrete probability measure.

Denote $k := |X| \geq 2$. Let $B_n$ be the set of all words in $X^*$ of length at most $n$. Thus $|B_n| = 1 + k + k^2 + \cdots + k^n = \frac{k^{n+1}-1}{k-1}$. Let $S \subseteq X^*$ be such that for each $w \in S$ the algorithm $M''$ terminates in time $f_1(|w|)$ and that $\lim_{n\to\infty} \frac{|S \cap B_n|}{|B_n|} = 1$ exponentially fast. Let $K = X^* - S$. Then there is $C > 0$ and $1 \leq q < k$ such that for every $n \geq 1$ we have $|K \cap B_n| \leq Cq^n$. For each $w \in X^*$ denote by $T(w)$ the time required for algorithm $M$ to reach a decision on input $w$. Thus for every $w \in S$ we have $T(w) \leq f_1(|w|)$ and for every $w \in X^*$ we have $T(w) \leq f(|w|)$.

Since $\mu$ is length-invariant, there is a function $d(n) \geq 0$ such that

$$\sum_{n=0}^{\infty} d(n) = 1$$

and such that for every $w \in X^*$ with $|w| = n$ we have $\mu(w) = \frac{d(n)}{k^n}$ where $k^n$ is the number of all words of length $n$ in $X^*$. In particular we have $d(n) \leq 1$ for each $n$.

We have:

$$\int_{X^*} \frac{T(w)}{f_1(|w|)}\mu(w) = \int_S \frac{T(w)}{f_1(|w|)}\mu(w) + \int_K \frac{T(w)}{f_1(|w|)}\mu(w).$$

Since for all $w \in S$ we have $T(w) \leq f_1(|w|)$, then

$$\int_S \frac{T(w)}{f_1(|w|)}\mu(w) \leq \int_S \frac{f_1(|w|)}{f_1(|w|)}\mu(w) = \mu(S) \leq 1 < \infty.$$

On the other hand for $w \in K$ with $|w| = n$ we have $T(w) \leq f(n)$ and $\mu(w) = d(n)\frac{1}{k^n}$. Hence $\frac{f(|w|)}{f_1(|w|)}d(|w|) \leq \frac{f(n)}{f_1(n)}$ where the function $\frac{f(n)}{f_1(n)}$ is subexponential. Let $K_n$ be the set of all $w \in K$ with $|w| = n$. Then $|K_n| \leq |K \cap B_n| \leq Cq^n$.

Hence

$$\int_K \frac{T(w)}{f_1(|w|)}\mu(w) = \sum_{n=0}^{\infty} \sum_{w \in K_n} \frac{T(w)}{f_1(|w|)}\mu(w) = \sum_{n=0}^{\infty} \sum_{w \in K_n} \frac{f(n)}{f_1(n)}d(n)\frac{1}{k^n} \leq$$

$$\leq C\sum_{n=0}^{\infty} \frac{f(n)}{f_1(n)}d(n)\frac{q^n}{k^n} \leq C\sum_{n=0}^{\infty} \frac{f(n)}{f_1(n)}\frac{q^n}{k^n} = C_0 < \infty$$

since $k > q$ and $\frac{f(n)}{f_1(n)}$ is subexponential. Therefore

$$\int_{X^*} \frac{T(w)}{f_1(|w|)}\mu(w) = \int_S \frac{T(w)}{f_1(|w|)}\mu(w) + \int_K \frac{T(w)}{f_1(|w|)}\mu(w) \leq 1 + C_0 < \infty.$$

Since $C_0$ does not depend on the choice of $\mu$, the statement of Proposition 3.2 follows. □

**Remark 3.3.** The proof of Proposition 3.2 is sufficiently robust to accommodate some other notions of average-case complexity. One of such definitions uses the Cauchy density $d(n) = \frac{1}{n^2}$ to define the notion of average-case complexity being at most polynomial of degree $m$. Namely, we say that, using the notations of Definition 1.4, *the algorithm $M$ solves $\mathcal{D}$ with Cauchy average-case time-complexity bounded by polynomial of degree $m$* if for any $\epsilon > 0$

$$\int_{X^*} \frac{T(w)}{|w|^{m-1+\epsilon}}\mu(w) = \sum_{w \in X^*} \frac{T(w)}{|w|^{m-1+\epsilon}}\mu(w) < \infty,$$

where $\mu(w) = \frac{1}{|w|^2 k^{|w|}}$.

The proof of Proposition 3.2 easily goes through to show that if $\mathcal{D}$ is decidable in subexponential time and strongly generically decidable in polynomial of degree $k$ time, then $\mathcal{D}$ is decidable with Cauchy average-case time-complexity bounded by polynomial of degree $k$.

Indeed, we have to analyze the same integral with $f_1(n) = n^{k-1+\epsilon}$. As before, we decompose this integral into two parts corresponding to $K$ and $S$ accordingly. For the $K$-part the functions $T(w)$ and $T(w)/|w|^{m-1+\epsilon}$ are subexponential and hence the $K$-integral is finite by the same argument as in the proof of Proposition 3.2. For the $S$-part we have $T(w) \le C|w|^m$ for any $w \in S$. Let $s_n$ be the number of words of length $n$ in $S$. Hence $s_n \le k^n$ and

$$\int_S \frac{T(w)}{|w|^{m-1+\epsilon}} \mu(w) \le C \int_S \frac{|w|^m}{|w|^{m-1+\epsilon}} \frac{1}{|w|^2 k^{|w|}} =$$
$$= C \sum_{n=1}^{\infty} \frac{1}{n^{1+\epsilon}} \frac{s_n}{k^n} \le C \sum_{n=1}^{\infty} \frac{1}{n^{1+\epsilon}} < \infty.$$

This implies that the whole integral over $X^*$ is finite, as required.

Proposition 3.2 together with the results of our previous paper [32] immediately implies the main results stated in the Introduction.

*Proof of Theorem A and Theorem C.* Suppose the assumptions of Theorem A are satisfied. By the results of [32] since $\bar{G}$ is non-amenable, for any finite generating set $A$ of $G$ the word problem $WP(G, A)$ for $G$ is solvable strongly generically with complexity $\mathcal{C}$. Together with Proposition 3.2 this implies the statement of Theorem A.

Any non-elementary hyperbolic group contains a free subgroup of rank two and thus is non-amenable. Moreover, by [3, 2, 30] the word problem in a hyperbolic group is solvable in linear time. Additionally, the word problem in an automatic group is solvable in quadratic time [21]. The statement of Corollary B now follows from Theorem A.

Suppose now that the assumptions of Theorem C hold. Since the Schreier graph of $\bar{G}$ over $\bar{H}$ is non-amenable, the result of [32] implies that for any finite generating set $A$ of $G$ the membership problem $MP(G, H, A)$ is solvable strongly generically with complexity $\mathcal{C}$. The statement of Theorem C now follows from Proposition 3.2.

Recall that the membership problem for a rational subgroup of a hyperbolic group is solvable in linear time. Also, the membership problem for a rational subgroup of an automatic group is solvable in quadratic time. Since non-elementary hyperbolic groups are non-amenable, Corollary D now follows from Theorem C. □

**Remark 3.4.** One may argue that the set of inputs for the word problem in a group $G = \langle A \rangle$ is the free group $F(A)$ rather than the set of all words (including those which are not freely reduced) in the alphabet $A \cup A^{-1}$. Then one would need to talk about average-case complexity of the word or the membership problem with respect to a length-invariant discrete probability measure $\mu : F(A) \to [0, 1]$. For an element $f \in F(A)$, $|f|_A$ is the length of the unique freely reduced word in $A \cup A^{-1}$ representing $f$. Here "length-invariant measure" would mean a measure such that for any $f_1, f_2 \in F(A)$ with $|f_1| = |f_2|$ we have $\mu(f_1) = \mu(f_2)$. Then the definition of average-case complexity would need to be modified so that the summation occurs over all elements of $F(A)$ rather than over all words in the language $(A \cup A^{-1})^*$

One can also define asymptotic density for a subset $S \subseteq F(A)$ similarly to Definition 2.1. The only difference is that the denominator in the fraction $\rho_n(S)$ would be the total number of elements $f \in F$ with $|f|_A \leq n$. All the results of [32] regarding the generic-case complexity of the word and the membership problem are proved in parallel for both approaches.

Thus all the results (together with exactly the same proofs) of this paper remain true in this alternative approach.

We give here a sample version of such an alternatively stated result:

**Theorem E.** *Let $G$ be a finitely presented group where the word problem is in deterministic time-complexity class $\mathcal{C}$. Suppose $G$ has a subgroup of finite index that possesses a non-amenable quotient group $\bar{G}$ where the word problem belongs to a time-complexity class $\mathcal{C} \subseteq \mathcal{S}ubExp$. Then for any finite generating set $A$ of $G$ the word problem for $G$ is solvable with average-case time-complexity in $\mathcal{C}$ uniformly relative to the family of all length-invariant discrete probability measures $\mu : F(A) \to [0, 1]$.*

## References

[1] S. I. Adian and V. G. Durnev, *Algorithmic problems for groups and semigroups,* Uspekhi Mat. Nauk **55** (2000), no. 2, 3–94; translation in Russian Math. Surveys **55** (2000), no. 2, 207–296

[2] J.Alonso, T.Brady, D.Cooper, V.Ferlini, M.Lustig, M.Mihalik, M.Shapiro and H.Short, *Notes on hyperbolic groups,* In: " Group theory from a geometrical viewpoint", Proceedings of the workshop held in Trieste, É. Ghys, A. Haefliger and A. Verjovsky (editors). World Scientific Publishing Co., 1991

[3] M. Anshel and B. Domanski, *The complexity of Dehn's algorithm for word problems in groups,* J. Algorithms **6** (1985), no. 4, 543–549

[4] G. Arzhantseva, *On groups in which subgroups with a fixed number of generators are free,*(Russian) Fundam. Prikl. Mat. **3** (1997), no. 3, 675–683

[5] G. Arzhantseva, *Generic properties of finitely presented groups and Howson's theorem,* Comm. Algebra **26** (1998), no. 11, 3783–3792

[6] G. Arzhantseva, *A property of subgroups of infinite index in a free group,* Proc. Amer. Math. Soc. **128** (2000), no. 11, 3205–3210

[7] G. Arzhantseva and A. Olshanskii, *Genericity of the class of groups in which subgroups with a lesser number of generators are free,* (Russian) Mat. Zametki **59** (1996), no. 4, 489–496

[8] L. Bartholdi, *Counting paths in graphs,* Enseign. Math. (2) **45** (1999), no. 1-2, 83–131

[9] B. Baumslag and S. J. Pride, *Groups with two more generators than relators,* J. London Math. Soc. (2) **17** (1978), no. 3, 425–426

[10] G. Baumslag, S. M. Gersten, M. Shapiro and H. Short, *Automatic groups and amalgams.* J. Pure Appl. Algebra **76** (1991), no. 3, 229–316

[11] G. Baumslag and C. F. Miller, III, *Experimenting and computing with infinite groups.* Groups and computation, II (New Brunswick, NJ, 1995), 19–30, DIMACS Ser. Discrete Math. Theoret. Comput. Sci., **28**, Amer. Math. Soc., Providence, RI, 1997.

[12] A. Borovik, A. G. Myasnikov, V. Shpilrain, *Measuring sets in infinite groups,* Contemp. Math., to appear.

[13] T. Ceccherini-Silberstein, R. Grigorchuck and P. de la Harpe, *Amenability and paradoxical decompositions for pseudogroups and discrete metric spaces,* (Russian) Tr. Mat. Inst. Steklova **224** (1999), Algebra. Topol. Differ. Uravn. i ikh Prilozh., 68–111; translation in Proc. Steklov Inst. Math., **224** (1999), no. 1, 57–97.

[14] C. Champetier, *Petite simplification dans les groupes hyperboliques*, Ann. Fac. Sci. Toulouse Math. (6) **3** (1994), no. 2, 161–221.

[15] C. Champetier, *Propriétés statistiques des groupes de présentation finie*, Adv. Math. **116** (1995), no. 2, 197–262.

[16] C. Champetier, *The space of finitely generated groups,* Topology **39** (2000), no. 4, 657–680

[17] P.-A. Cherix and A. Valette, *On spectra of simple random walks on one-relator groups,* With an appendix by Paul Jolissaint. Pacific J. Math. **175** (1996), no. 2, 417–438

[18] P.-A. Cherix and G. Schaeffer, *An asymptotic Freiheitssatz for finitely generated groups,* Enseign. Math. (2) **44** (1998), no. 1-2, 9–22

[19] J. M. Cohen, *Cogrowth and amenability of discrete groups*, J. Funct. Anal. **48** (1982), 301–309.

[20] M. Coornaert, T. Delzant, and A. Papadopoulos, *Géométrie et théorie des groupes. Les groupes hyperboliques de Gromov.* Lecture Notes in Mathematics, 1441; Springer-Verlag, Berlin, 1990

[21] D. Epstein, J. Cannon, D. Holt, S. Levy, M. Paterson, W. Thurston, *Word Processing in Groups,* Jones and Bartlett, Boston, 1992

[22] B. Farb, *Automatic groups: a guided tour,* Enseign. Math. (2) **38** (1992), no. 3-4, 291–313

[23] S. Gersten, *Introduction to hyperbolic and automatic groups.* Summer School in Group Theory in Banff, 1996, 45–70, CRM Proc. Lecture Notes, **17**, Amer. Math. Soc., Providence, RI, 1999

[24] S. Gersten and H. Short, *Rational subgroups of bi-automatic groups*

[25] E. Ghys and P. de la Harpe (editors), *Sur les groupes hyperboliques d'aprés Mikhael Gromov*, Birkhäuser, Progress in Mathematics series, vol. **83**, 1990.

[26] R. Grigorchuk, *Symmetrical random walks on discrete groups,* Multicomponent random systems, pp. 285–325, Adv. Probab. Related Topics, **6**, Dekker, New York, 1980

[27] M. Gromov, *Hyperbolic groups*, Essays in group theory, Springer, New York, 1987, pp. 75–263.

[28] M. Gromov, *Random walks in random groups*, preprint, 2001

[29] ———, *Asymptotic invariants of infinite groups*, Geometric group theory, Vol. 2 (Sussex, 1991), Cambridge Univ. Press, Cambridge, 1993, pp. 1–295.

[30] D. Holt, *Word-hyperbolic groups have real-time word problem,* Internat. J. Algebra Comput. *10* (2000), no. 2, 221–227

[31] I. Kapovich, *The non-amenability of Schreier graphs for infinite index quasiconvex subgroups of hyperbolic groups*, to appear in Ensign. Math.

[32] I. Kapovich, A. Myasnikov, V. Shpilrain and P. E. Schupp, *Generic-case complexity and decision problems in group theory and random walks*, preprint, 2002

[33] L. Levin, *Average case complete problems,* SIAM J. Comput. **15** (1986), 285–286

[34] R. Lipton and Y. Zalcstein, *Word problems solvable in logspace,* J. Assoc. Comput. Mach. **24** (1977), no. 3, 522–526

[35] R. C. Lyndon and P. E. Schupp, *Combinatorial Group Theory*, Ergebnisse der Mathematik, band 89, Springer 1977. Reprinted in the Springer Classics in Mathematics series, 2000.

[36] C. F. Miller III, *Decision problems for groups- Survey and reflections*. in *Algorithms and Classification in Combinatorial Group Theory*, G. Bamuslag and C.F. Miller III, editors, (1992), Springer, 1–60.

[37] C. Papadimitriou, *Computation Complexity*, (1994), Addison-Wesley, Reading.
[38] A. Yu. Ol'shanskii, *Almost every group is hyperbolic*, Internat. J. Algebra Comput. **2** (1992), no. 1, 1–17.
[39] H. Short, *An introduction to automatic groups.* Semigroups, formal languages and groups (York, 1993), 233–253, NATO Adv. Sci. Inst. Ser. C Math. Phys. Sci., **466**, Kluwer Acad. Publ., Dordrecht, 1995
[40] A. M. Vershik, *Dynamic theory of growth in groups: entropy, boundaries, examples* Uspekhi Mat. Nauk **55** (2000), no. 4(334), 59–128; translation in Russian Math. Surveys **55** (2000), no. 4, 667–733
[41] J. Wang, *Average-case completeness of a word problem in groups,* Proc. of the 27-th Annual Symposium on Theory of Computing, ACM Press, New York, 1995, 325–334
[42] J. Wang, *Average-case computational complexity theory,* Complexity Theory Retrospective, II. Springer-Verlag, New York, 1997, 295–334
[43] J. Wang, *Distributional word problem for groups,* SIAM J. Comput. **28** (1999), no. 4, 1264–1283
[44] W. Woess, *Cogrowth of groups and simple random walks*, Arch. Math. **41** (1983), 363–370.
[45] W. Woess, *Random walks on infinite graphs and groups - a survey on selected topics,* Bull. London Math. Soc. **26** (1994), 1–60.

Department of Mathematics, University of Illinois at Urbana-Champaign, 1409 West Green Street, Urbana, IL 61801, USA
  *E-mail address*: kapovich@math.uiuc.edu

Department of Mathematics, The City College of New York, New York, NY 10031
  *E-mail address*: alexeim@att.net

Department of Mathematics, University of Illinois at Urbana-Champaign, 1409 West Green Street, Urbana, IL 61801, USA
  *E-mail address*: schupp@math.uiuc.edu

Department of Mathematics, The City College of New York, New York, NY 10031
  *E-mail address*: shpil@groups.sci.ccny.cuny.edu